# ASYMPTOTIC BEHAVIOR OF WEIGHTED QUADRATIC VARIATIONS OF FRACTIONAL BROWNIAN MOTION: THE CRITICAL CASE $H = 1/4$

By Ivan Nourdin and Anthony Réveillac

*Université Paris VI and Humboldt-Universität zu Berlin*

We derive the asymptotic behavior of weighted quadratic variations of fractional Brownian motion $B$ with Hurst index $H = 1/4$. This completes the only missing case in a very recent work by I. Nourdin, D. Nualart and C. A. Tudor. Moreover, as an application, we solve a recent conjecture of K. Burdzy and J. Swanson on the asymptotic behavior of the Riemann sums with alternating signs associated to $B$.

**1. Introduction.** Let $B^H$ be a fractional Brownian motion with the Hurst index $H \in (0, 1)$. Drawing on the seminal works by Breuer and Major [1], Dobrushin and Major [5], Giraitis and Surgailis [6] or Taqqu [24], it is well known that:

- if $H \in (0, \frac{3}{4})$, then

(1.1) $$\frac{1}{\sqrt{n}} \sum_{k=0}^{n-1} [n^{2H}(B^H_{(k+1)/n} - B^H_{k/n})^2 - 1] \xrightarrow[n \to \infty]{\text{Law}} \mathcal{N}(0, C_H^2);$$

- if $H = \frac{3}{4}$, then

(1.2) $$\frac{1}{\sqrt{n \log n}} \sum_{k=0}^{n-1} [n^{3/2}(B^{3/4}_{(k+1)/n} - B^{3/4}_{k/n})^2 - 1] \xrightarrow[n \to \infty]{\text{Law}} \mathcal{N}(0, C_{3/4}^2);$$

- if $H \in (\frac{3}{4}, 1)$, then

(1.3) $$n^{1-2H} \sum_{k=0}^{n-1} [n^{2H}(B^H_{(k+1)/n} - B^H_{k/n})^2 - 1] \xrightarrow[n \to \infty]{\text{Law}} \text{``Rosenblatt r.v.''}$$

Here, $C_H > 0$ denotes a constant depending only on $H$ which can be computed explicitly. Moreover, the term "Rosenblatt r.v." denotes a random variable whose distribution is the same as that of the Rosenblatt process $Z$ at time one [see (1.9) below].

Now, let $f$ be a (regular enough) real function. Very recently, the asymptotic behavior of

$$\sum_{k=0}^{n-1} f(B_{k/n}^H)[n^{2H}(B_{(k+1)/n}^H - B_{k/n}^H)^2 - 1] \quad (1.4)$$

received a lot of attention (see [7, 13, 14, 15, 17]) (see also the related works [18, 21, 22, 23]). The initial motivation of such a study was to derive the exact rates of convergence of some approximation schemes associated with scalar stochastic differential equations driven by $B^H$ (see [7, 13, 14] for precise statements). But it turned out that it was also interesting because it highlighted new phenomena with respect to (1.1), (1.2) and (1.3). Indeed, in the study of the asymptotic behavior of (1.4), a new critical value ($H = \frac{1}{4}$) appears. More precisely:

- if $H < \frac{1}{4}$, then

$$(1.5) \quad n^{2H-1} \sum_{k=0}^{n-1} f(B_{k/n}^H)[n^{2H}(B_{(k+1)/n}^H - B_{k/n}^H)^2 - 1]$$

$$\xrightarrow[n \to \infty]{L^2} \frac{1}{4} \int_0^1 f''(B_s^H) \, ds;$$

- if $\frac{1}{4} < H < \frac{3}{4}$, then

$$(1.6) \quad \frac{1}{\sqrt{n}} \sum_{k=0}^{n-1} f(B_{k/n}^H)[n^{2H}(B_{(k+1)/n}^H - B_{k/n}^H)^2 - 1]$$

$$\xrightarrow[n \to \infty]{\text{Law}} C_H \int_0^1 f(B_s^H) \, dW_s$$

for $W$ a standard Brownian motion independent of $B^H$;
- if $H = \frac{3}{4}$, then

$$(1.7) \quad \frac{1}{\sqrt{n \log n}} \sum_{k=0}^{n-1} f(B_{k/n}^{3/4})[n^{3/2}(B_{(k+1)/n}^{3/4} - B_{k/n}^{3/4})^2 - 1]$$

$$\xrightarrow[n \to \infty]{\text{Law}} C_{3/4} \int_0^1 f(B_s^{3/4}) \, dW_s$$

for $W$ a standard Brownian motion independent of $B^{3/4}$;

- if $H > \frac{3}{4}$ then

$$(1.8) \quad n^{1-2H} \sum_{k=0}^{n-1} f(B_{k/n}^H)[n^{2H}(B_{(k+1)/n}^H - B_{k/n}^H)^2 - 1] \xrightarrow[n\to\infty]{L^2} \int_0^1 f(B_s^H) \, dZ_s$$

for $Z$ the Rosenblatt process defined by

$$(1.9) \quad Z_s = I_2^X(L_s),$$

where $I_2^X$ denotes the double stochastic integral with respect to the Wiener process $X$ given by the transfer equation (2.3) and where, for every $s \in [0,1]$, $L_s$ is the symmetric square integrable kernel given by

$$L_s(y_1, y_2) = \frac{1}{2} \mathbf{1}_{[0,s]^2}(y_1, y_2) \int_{y_1 \vee y_2}^s \frac{\partial K_H}{\partial u}(u, y_1) \frac{\partial K_H}{\partial u}(u, y_2) \, du.$$

Even if it is not completely obvious at first glance, convergences (1.1) and (1.5) agree. Indeed, since $2H - 1 < -\frac{1}{2}$ if and only if $H < \frac{1}{4}$, (1.5) is actually a particular case of (1.1) when $f \equiv 1$. The convergence (1.5) is proved in [15] while cases (1.6), (1.7) and (1.8) are proved in [17]. On the other hand, notice that the relation (1.5) to (1.8) do not cover the critical case $H = \frac{1}{4}$. Our first main result (see below) completes this important missing case.

THEOREM 1.1. *If $H = \frac{1}{4}$, then*

$$
\begin{aligned}
&\frac{1}{\sqrt{n}} \sum_{k=0}^{n-1} f(B_{k/n}^{1/4})[\sqrt{n}(B_{(k+1)/n}^{1/4} - B_{k/n}^{1/4})^2 - 1] \\
&\xrightarrow[n\to\infty]{\text{Law}} C_{1/4} \int_0^1 f(B_s^{1/4}) \, dW_s + \frac{1}{4} \int_0^1 f''(B_s^{1/4}) \, ds
\end{aligned}
\tag{1.10}
$$

*for $W$ a standard Brownian motion independent of $B^{1/4}$ and where*

$$C_{1/4} = \sqrt{\frac{1}{2} \sum_{p=-\infty}^{\infty} (\sqrt{|p+1|} + \sqrt{|p-1|} - 2\sqrt{|p|})^2} \approx 1535.$$

Here, it is interesting to compare the obtained limit in (1.10) with those obtained in the recent work [18]. In [18], the authors study the asymptotic behavior of (1.4) when the fractional Brownian motion $B^H$ is replaced by an *iterated Brownian motion* $Z$, that is, the process defined by $Z_t = X(Y_t)$, $t \in [0,1]$, with $X$ and $Y$, two independent Brownian motions. Iterated Brownian motion $Z$ is self similar of index $H = 1/4$ and has stationary increments. Hence although it is not Gaussian, $Z$ is "close" to the fractional Brownian motion $B^{1/4}$. For $Z$ instead of $B^{1/4}$, it is proved in [18] that the correctly renormalized weighted quadratic variation [which is not exactly defined as in (1.4), but rather by means of a *random* partition composed of Brownian hitting times] converges in law toward the *weighted Brownian motion in random scenery* at time one, defined as

$$\sqrt{2}\int_{-\infty}^{+\infty} f(X_x)L_1^x(Y)\,dW_x,$$

compare with the right-hand side of (1.10). Here, $\{L_t^x(Y)\}_{x\in\mathbb{R},t\in[0,1]}$ stands for the jointly continuous version of the local time process of $Y$, while $W$ denotes a two-sided standard Brownian motion independent of $X$ and $Y$.

From now on, we will only work with a fractional Brownian motion of the Hurst index $H=\frac{1}{4}$. This particular value of $H$ is important because the fractional Brownian motion with the Hurst index $H=\frac{1}{4}$ has a remarkable physical interpretation in terms of particle systems. Indeed, if one considers an infinite number of particles, initially placed on the real line according to a Poisson distribution, performing independent Brownian motions and undergoing "elastic" collisions, then the trajectory of a fixed particle (after rescaling) converges to a fractional Brownian motion with the Hurst index $H=\frac{1}{4}$. This striking fact has been first pointed out by Harris in [11], and has been rigorously proven in [4] (see also references therein).

Now let us explain an interesting consequence of a slight modification of Theorem 1.1 toward the first step in a construction of a stochastic calculus with respect to $B^{1/4}$. As it is nicely explained by Swanson in [23], there are (at least) two kinds of Stratonovitch-type Riemann sums that one can consider in order to define $\int_0^1 f(B_s^{1/4})\circ dB_s^{1/4}$ when $f$ is a real (smooth enough) function. The first one corresponds to the so-called "trapezoid rule" and is given by

$$S_n(f)=\sum_{k=0}^{n-1}\frac{f(B_{k/n}^{1/4})+f(B_{(k+1)/n}^{1/4})}{2}(B_{(k+1)/n}^{1/4}-B_{k/n}^{1/4}).$$

The second one corresponds to the so-called "midpoint rule" and is given by

$$T_n(f)=\sum_{k=1}^{\lfloor n/2\rfloor} f(B_{(2k-1)/n}^{1/4})(B_{(2k)/n}^{1/4}-B_{(2k-2)/n}^{1/4}).$$

By Theorem 3 in [17] (see also [3, 8, 9]), we have that

$$\int_0^1 f'(B_s^{1/4})\,d^\circ B_s^{1/4}:=\lim_{n\to\infty} S_n(f') \qquad \text{exists in probability}$$

and verifies the following *classical* change of variable formula:

(1.11)  $$\int_0^1 f'(B_s^{1/4})\,d^\circ B_s^{1/4}=f(B_1^{1/4})-f(0).$$

On the other hand, it is quoted in [23] that Burdzy and Swanson *conjectured*[1] that

$$\int_0^1 f'(B_s^{1/4})\, d^\star B_s^{1/4} := \lim_{n\to\infty} T_n(f') \qquad \text{exists in law}$$

and verifies, this time, the following *nonclassical* change of variable formula:

$$(1.13) \qquad \int_0^1 f'(B_s^{1/4})\, d^\star B_s^{1/4} \stackrel{\text{Law}}{=} f(B_1^{1/4}) - f(0) - \frac{\kappa}{2}\int_0^1 f''(B_s^{1/4})\, dW_s,$$

where $\kappa$ is an explicit universal constant, and $W$ denotes a standard Brownian motion independent of $B^{1/4}$. Our second main result is the following:

THEOREM 1.2. *The conjecture of Burdzy and Swanson is true. More precisely, (1.13) holds for any real function $f:\mathbb{R}\to\mathbb{R}$ verifying* ($H_9$) *(see Section 3 below).*

The rest of the paper is organized as follows. In Section 2, we recall some notions concerning fractional Brownian motion. In Section 3, we prove Theorem 1.1. In Section 4, we prove Theorem 1.2.

*Note*: Just after this paper was put on the ArXiv, Burdzy and Swanson informed us that, in their manuscript [2], prepared independently and at the same time as ours, they also proved Theorem 1.2 by using a completely different route.

**2. Preliminaries and notation.** We begin by briefly recalling some basic facts about stochastic calculus with respect to a fractional Brownian motion. We refer to [19, 20] for further details. Let $B^H = (B_t^H)_{t\in[0,1]}$ be a fractional

---

[1]Actually, Burdzy and Swanson have conjectured (1.13) not for the fractional Brownian motion $B^{1/4}$ but for the process $F$ defined by

$$(1.12) \qquad F_t = u(t,0), \qquad t \in [0,1],$$

where

$$u_t = \tfrac{1}{2}u_{xx} + \dot{W}(t,x), t \in [0,1], x \in \mathbb{R} \qquad \text{with initial condition } u(0,x) = 0, x \in \mathbb{R}.$$

(Here, as usual, $\dot{W}$ denotes the space-time white noise on $[0,1] \times \mathbb{R}$.) It is immediately checked that $F$ is a centered Gaussian process with covariance function

$$E(F_s F_t) = \frac{1}{\sqrt{2\pi}}(\sqrt{t+s} - \sqrt{|t-s|}),$$

so that $F$ is a *bifractional Brownian motion* of indices $\frac{1}{2}$ and $\frac{1}{2}$ in the sense of Houdré and Villa [10]. Using the main result of [12], we have that $B^{1/4}$ and $F$ actually differ only from a process with absolutely continuous trajectories. As a direct consequence, using a Girsanov-type transformation, it is equivalent to prove (1.13) either for $B^{1/4}$ or for $F$.

Brownian motion with the Hurst parameter $H \in (0, \frac{1}{2})$. That is, $B^H$ is a centered Gaussian process with the covariance function

$$R_H(s,t) = \tfrac{1}{2}(t^{2H} + s^{2H} - |t-s|^{2H}). \tag{2.1}$$

We denote by $\mathscr{E}$ the set of step $\mathbb{R}$-valued functions on $[0,1]$. Let $\mathfrak{H}$ be the Hilbert space defined as the closure of $\mathscr{E}$ with respect to the scalar product

$$\langle \mathbf{1}_{[0,t]}, \mathbf{1}_{[0,s]} \rangle_{\mathfrak{H}} = R_H(t,s).$$

The covariance kernel $R_H(t,s)$ introduced in (2.1) can be written as

$$R_H(t,s) = \int_0^{s \wedge t} K_H(s,u) K_H(t,u)\, du,$$

where $K_H(t,s)$ is the square integrable kernel defined, for $0 < s < t$, by

$$K_H(t,s) = c_H \bigg[ \bigg(\frac{t}{s}\bigg)^{H-1/2} (t-s)^{H-1/2} \\ - (H-1/2) s^{1/2-H} \int_s^t u^{H-3/2}(u-s)^{H-1/2}\, du \bigg], \tag{2.2}$$

where $c_H^2 = 2H(1-2H)^{-1}\beta(1-2H, H+1/2)^{-1}$ ($\beta$ denotes the Beta function). By convention, we set $K_H(t,s) = 0$ if $s \geq t$.

Let $\mathcal{K}_H^* : \mathscr{E} \to L^2([0,1])$ be the linear operator defined by

$$\mathcal{K}_H^*(\mathbf{1}_{[0,t]}) = K_H(t, \cdot).$$

The following equality holds for any $s, t \in [0, 1]$:

$$\langle \mathbf{1}_{[0,t]}, \mathbf{1}_{[0,s]} \rangle_{\mathfrak{H}} = \langle \mathcal{K}_H^* \mathbf{1}_{[0,t]}, \mathcal{K}_H^* \mathbf{1}_{[0,s]} \rangle_{L^2([0,1])} = E(B_t^H B_s^H);$$

hence $\mathcal{K}_H^*$ provides an isometry between the Hilbert spaces $\mathfrak{H}$ and a closed subspace of $L^2([0,1])$. The process $X = (X_t)_{t \in [0,1]}$ defined by

$$X_t = B^H((\mathcal{K}_H^*)^{-1}(\mathbf{1}_{[0,t]})) \tag{2.3}$$

is a standard Brownian motion, and the process $B^H$ has an integral representation of the form

$$B_t^H = \int_0^t K_H(t,s)\, dX_s.$$

Let $\mathscr{S}$ be the set of all smooth cylindrical random variables, that is, of the form

$$F = \psi(B_{t_1}^H, \ldots, B_{t_m}^H), \tag{2.4}$$

where $m \geq 1$, $\psi : \mathbb{R}^m \to \mathbb{R} \in \mathscr{C}_b^\infty$ and $0 \leq t_1 < \cdots < t_m \leq 1$. The Malliavin derivative of $F$ with respect to $B^H$ is the element of $L^2(\Omega, \mathfrak{H})$ defined by

$$D_s F = \sum_{i=1}^m \frac{\partial \psi}{\partial x_i}(B_{t_1}^H, \ldots, B_{t_m}^H) \mathbf{1}_{[0,t_i]}(s), \qquad s \in [0,1].$$

In particular $D_s B_t^H = \mathbf{1}_{[0,t]}(s)$. For any integer $k \geq 1$, we denote by $\mathbb{D}^{k,2}$ the closure of the set of smooth random variables with respect to the norm

$$\|F\|_{k,2}^2 = \mathrm{E}[F^2] + \sum_{j=1}^k \mathrm{E}[|D^j F|^2_{\mathfrak{H}^{\otimes j}}].$$

The Malliavin derivative $D$ verifies the chain rule: if $\varphi : \mathbb{R}^n \to \mathbb{R}$ is $\mathscr{C}_b^1$, and if $(F_i)_{i=1,\ldots,n}$ is a sequence of elements of $\mathbb{D}^{1,2}$, then $\varphi(F_1,\ldots,F_n) \in \mathbb{D}^{1,2}$, and we have, for any $s \in [0,1]$,

$$D_s \varphi(F_1,\ldots,F_n) = \sum_{i=1}^n \frac{\partial \varphi}{\partial x_i}(F_1,\ldots,F_n) D_s F_i.$$

The divergence operator $I$ is the adjoint of the derivative operator $D$. If a random variable $u \in \mathrm{L}^2(\Omega, \mathfrak{H})$ belongs to the domain of the divergence operator, that is, if it verifies

$$|\mathrm{E}\langle DF, u\rangle_{\mathfrak{H}}| \leq c_u \|F\|_{\mathrm{L}^2} \qquad \text{for any } F \in \mathscr{S},$$

then $I(u)$ is defined by the duality relationship

$$\mathrm{E}(FI(u)) = \mathrm{E}(\langle DF, u\rangle_{\mathfrak{H}})$$

for every $F \in \mathbb{D}^{1,2}$.

For every $n \geq 1$, let $\mathcal{H}_n$ be the $n$th Wiener chaos of $B^H$, that is, the closed linear subspace of $L^2(\Omega)$ generated by the random variables $\{H_n(B^H(h)), h \in H, |h|_{\mathfrak{H}} = 1\}$ where $H_n$ is the $n$th Hermite polynomial. The mapping $I_n(h^{\otimes n}) = n! H_n(B^H(h))$ provides a linear isometry between the symmetric tensor product $\mathfrak{H}^{\odot n}$ and $\mathcal{H}_n$. For $H = \frac{1}{2}$, $I_n$ coincides with the multiple stochastic integral. The following duality formula holds

$$(2.5) \qquad E(FI_n(h)) = E(\langle D^n F, h\rangle_{\mathfrak{H}^{\otimes n}})$$

for any element $h \in \mathfrak{H}^{\odot n}$ and any random variable $F \in \mathbb{D}^{n,2}$. Let $\{e_k, k \geq 1\}$ be a complete orthonormal system in $\mathfrak{H}$. Given $f \in \mathfrak{H}^{\odot p}$ and $g \in \mathfrak{H}^{\odot q}$, for every $r = 0, \ldots, p \wedge q$, the $r$th contraction of $f$ and $g$ is the element of $\mathfrak{H}^{\otimes(p+q-2r)}$ defined as

$$f \otimes_r g = \sum_{i_1,\ldots,i_r=1}^\infty \langle f, e_{i_1} \otimes \cdots \otimes e_{i_r}\rangle_{\mathfrak{H}^{\otimes r}} \otimes \langle g, e_{i_1} \otimes \cdots \otimes e_{i_r}\rangle_{\mathfrak{H}^{\otimes r}}.$$

Note that $f \otimes_0 g = f \otimes g$ equals the tensor product of $f$ and $g$ while for $p = q$, $f \otimes_p g = \langle f, g\rangle_{\mathfrak{H}^{\otimes p}}$. Finally, we mention the useful following multiplication formula: if $f \in \mathfrak{H}^{\odot p}$ and $g \in \mathfrak{H}^{\odot q}$, then

$$(2.6) \qquad I_p(f) I_q(g) = \sum_{r=0}^{p \wedge q} r! \binom{p}{r} \binom{q}{r} I_{p+q-2r}(f \otimes_r g).$$

**3. Proof of Theorem 1.1.** In this section, $B = B^{1/4}$ denotes a fractional Brownian motion with the Hurst index $H = 1/4$. Let

$$G_n := \frac{1}{\sqrt{n}} \sum_{k=0}^{n-1} f(B_{k/n})[\sqrt{n}(B_{(k+1)/n} - B_{k/n})^2 - 1], \qquad n \geq 1.$$

For $k = 0, \ldots, n-1$ and $t \in [0,1]$, we set

$$\delta_{k/n} := \mathbf{1}_{[k/n,(k+1)/n]} \quad \text{and} \quad \varepsilon_t := \mathbf{1}_{[0,t]}.$$

The relations between Hermite polynomials and multiple stochastic integrals (see Section 2) allow one to write

$$\sqrt{n}(B_{(k+1)/n} - B_{k/n})^2 - 1 = \sqrt{n} I_2(\delta_{k/n}^{\otimes 2}).$$

As a consequence,

$$G_n = \sum_{k=0}^{n-1} f(B_{k/n}) I_2(\delta_{k/n}^{\otimes 2}).$$

In the sequel, we will need the following assumption:

*Hypothesis* $(\mathrm{H}_q)$. The function $f : \mathbb{R} \to \mathbb{R}$ belongs to $\mathscr{C}^q$ and is such that

$$\sup_{t \in [0,1]} E(|f^{(i)}(B_t)|^p) < \infty$$

for any $p \geq 1$ and $i \in \{0, \ldots, q\}$.

We begin by the following technical lemma:

LEMMA 3.1. *Let $n \geq 1$ and $k = 0, \ldots, n-1$. We have:*

(i) $|E(B_r(B_t - B_s))| \leq \sqrt{t-s}$ *for any $r \in [0,1]$ and $0 \leq s \leq t \leq 1$;*
(ii) $\sup_{t \in [0,1]} \sum_{k=0}^{n-1} |\langle \varepsilon_t; \delta_{k/n} \rangle_{\mathfrak{H}}| \underset{n \to \infty}{=} O(1)$;
(iii) $\sum_{k,j=0}^{n-1} |\langle \varepsilon_{j/n}, \delta_{k/n} \rangle_{\mathfrak{H}}| \underset{n \to \infty}{=} O(n)$;
(iv) $|\langle \varepsilon_{k/n}, \delta_{k/n} \rangle_{\mathfrak{H}}^2 - \frac{1}{4n}| \leq \frac{\sqrt{k+1} - \sqrt{k}}{2n}$; *consequently* $\sum_{k=0}^{n-1} |\langle \varepsilon_{k/n}, \delta_{k/n} \rangle_{\mathfrak{H}}^2 - \frac{1}{4n}| \underset{n \to \infty}{\longrightarrow} 0$.

PROOF.

(i) We have

$$E(B_r(B_t - B_s)) = \tfrac{1}{2}(\sqrt{t} - \sqrt{s}) + \tfrac{1}{2}(\sqrt{|s-r|} - \sqrt{|t-r|}).$$

Using the classical inequality $|\sqrt{|b|} - \sqrt{|a|}| \leq \sqrt{|b-a|}$, the desired result follows.

(ii) Observe that
$$\langle \varepsilon_t, \delta_{k/n} \rangle_{\mathfrak{H}} = \frac{1}{2\sqrt{n}}(\sqrt{k+1} - \sqrt{k} - \sqrt{|k+1-nt|} + \sqrt{|k-nt|}).$$

Consequently, we have
$$\sum_{k=0}^{n-1} |\langle \varepsilon_t, \delta_{k/n} \rangle_{\mathfrak{H}}|$$
$$\leq \frac{1}{2} + \frac{1}{2\sqrt{n}} \Bigg( \sum_{k=0}^{\lfloor nt \rfloor - 1} \sqrt{nt-k} - \sqrt{nt-k-1}$$
$$+ \sqrt{\lfloor nt \rfloor + 1 - nt} - \sqrt{nt - \lfloor nt \rfloor}$$
$$+ \sum_{k=\lfloor nt \rfloor + 1}^{n-1} \sqrt{k-nt} - \sqrt{k+1-nt} \Bigg).$$

The desired conclusion follows easily.

(iii) It is a direct consequence of (ii),
$$\sum_{k,j=0}^{n-1} |\langle \varepsilon_{j/n}, \delta_{k/n} \rangle_{\mathfrak{H}}| \leq n \sup_{j=0,\ldots,n-1} \sum_{k=0}^{n-1} |\langle \varepsilon_{j/n}, \delta_{k/n} \rangle_{\mathfrak{H}}|$$
$$\underset{n \to \infty}{=} O(n).$$

(iv) We have
$$\left| \langle \varepsilon_{k/n}, \delta_{k/n} \rangle_{\mathfrak{H}}^2 - \frac{1}{4n} \right|$$
$$= \frac{1}{4n}(\sqrt{k+1} - \sqrt{k})|\sqrt{k+1} - \sqrt{k} - 2|.$$

Thus the desired bound is immediately checked by using $0 \leq \sqrt{x+1} - \sqrt{x} \leq 1$ available for $x \geq 0$. □

The main result of this section is the following:

THEOREM 3.2. *Under Hypothesis* $(H_4)$, *we have*
$$G_n \xrightarrow[n \to \infty]{\text{Law}} C_{1/4} \int_0^1 f(B_s) \, dW_s + \frac{1}{4} \int_0^1 f''(B_s) \, ds,$$
*where* $W = (W_t)_{t \in [0,1]}$ *is a standard Brownian motion independent of $B$ and*
$$C_{1/4} := \sqrt{\frac{1}{2} \sum_{p=-\infty}^{\infty} (\sqrt{|p+1|} + \sqrt{|p-1|} - 2\sqrt{|p|})^2} \approx 1535.$$

PROOF. This proof is mainly inspired by the first draft of [16]. Throughout the proof, $C$ will denote a constant depending only on $\|f^{(a)}\|_\infty$, $a = 0, 1, 2, 3, 4$, which can differ from one line to another.

STEP 1. We begin the proof by showing the following limits:

$$\lim_{n\to\infty} E(G_n) = \frac{1}{4}\int_0^1 E(f''(B_s))\,ds \tag{3.1}$$

and

$$\lim_{n\to\infty} E(G_n^2) = C_{1/4}^2 \int_0^1 E(f^2(B_s))\,ds + \frac{1}{16}E\left(\int_0^1 f''(B_s)\,ds\right)^2. \tag{3.2}$$

PROOF OF (3.1). We can write

$$\begin{aligned}
E(G_n) &= \sum_{k=0}^{n-1} E(f(B_{k/n})I_2(\delta_{k/n}^{\otimes 2})) \\
&= \sum_{k=0}^{n-1} E(\langle D^2(f(B_{k/n})), \delta_{k/n}^{\otimes 2}\rangle_{\mathfrak{H}}) \\
&= \sum_{k=0}^{n-1} E(f''(B_{k/n}))\langle \varepsilon_{k/n}, \delta_{k/n}\rangle_{\mathfrak{H}}^2 \\
&= \frac{1}{4n}\sum_{k=0}^{n-1} E(f''(B_{k/n})) \\
&\quad + \sum_{k=0}^{n-1} E(f''(B_{k/n}))\left(\langle \varepsilon_{k/n}, \delta_{k/n}\rangle_{\mathfrak{H}}^2 - \frac{1}{4n}\right) \\
&\xrightarrow[n\to\infty]{} \frac{1}{4}\int_0^1 E(f''(B_s))\,ds \qquad \text{by Lemma 3.1(iv) and under (H}_4\text{).}
\end{aligned}$$

PROOF OF (3.2). By the multiplication formula (2.6), we have

$$I_2(\delta_{j/n}^{\otimes 2})I_2(\delta_{k/n}^{\otimes 2}) = I_4(\delta_{j/n}^{\otimes 2} \otimes \delta_{k/n}^{\otimes 2}) + 4I_2(\delta_{j/n} \otimes \delta_{k/n})\langle \delta_{j/n}, \delta_{k/n}\rangle_{\mathfrak{H}} + 2\langle \delta_{j/n}, \delta_{k/n}\rangle_{\mathfrak{H}}^2. \tag{3.3}$$

Thus

$$E(G_n^2) = \sum_{j,k=0}^{n-1} E(f(B_{j/n})f(B_{k/n})I_2(\delta_{j/n}^{\otimes 2})I_2(\delta_{k/n}^{\otimes 2}))$$

$$= \sum_{j,k=0}^{n-1} E(f(B_{j/n})f(B_{k/n})I_4(\delta_{j/n}^{\otimes 2} \otimes \delta_{k/n}^{\otimes 2}))$$

$$+ 4 \sum_{j,k=0}^{n-1} E(f(B_{j/n})f(B_{k/n})I_2(\delta_{j/n} \otimes \delta_{k/n}))\langle \delta_{j/n}, \delta_{k/n}\rangle_{\mathfrak{H}}$$

$$+ 2 \sum_{j,k=0}^{n-1} E(f(B_{j/n})f(B_{k/n}))\langle \delta_{j/n}, \delta_{k/n}\rangle_{\mathfrak{H}}^2$$

$$= A_n + B_n + C_n.$$

Using the Malliavin integration by parts formula (2.5), $A_n$ can be expressed as follows:

$$A_n = \sum_{j,k=0}^{n-1} E(\langle D^4(f(B_{j/n})f(B_{k/n})), \delta_{j/n}^{\otimes 2} \otimes \delta_{k/n}^{\otimes 2}\rangle_{\mathfrak{H}^{\otimes 4}})$$

$$= 24 \sum_{j,k=0}^{n-1} \sum_{a+b=4} E(f^{(a)}(B_{j/n})f^{(b)}(B_{k/n}))\langle \varepsilon_{j/n}^{\otimes a}\widetilde{\otimes}\varepsilon_{k/n}^{\otimes b}, \delta_{j/n}^{\otimes 2} \otimes \delta_{k/n}^{\otimes 2}\rangle_{\mathfrak{H}^{\otimes 4}}.$$

In the previous sum, each term is negligible except

$$\sum_{j,k=0}^{n-1} E(f''(B_{j/n})f''(B_{k/n}))\langle \varepsilon_{j/n}, \delta_{j/n}\rangle_{\mathfrak{H}}^2 \langle \varepsilon_{k/n}, \delta_{k/n}\rangle_{\mathfrak{H}}^2$$

$$= E\left(\left[\sum_{k=0}^{n-1} f''(B_{k/n})\langle \varepsilon_{k/n}, \delta_{k/n}\rangle_{\mathfrak{H}}^2\right]^2\right)$$

$$= E\left(\left[\frac{1}{4n}\sum_{k=0}^{n-1} f''(B_{k/n}) + \sum_{k=0}^{n-1} f''(B_{k/n})\left(\langle \varepsilon_{k/n}, \delta_{k/n}\rangle_{\mathfrak{H}}^2 - \frac{1}{4n}\right)\right]^2\right)$$

$$\xrightarrow[n\to\infty]{} E\left(\left[\frac{1}{4}\int_0^1 f''(B_s)\,ds\right]^2\right) \quad \text{by Lemma 3.1(iv) and under (H}_4\text{).}$$

The other terms appearing in $A_n$ make no contribution to the limit. Indeed, they have the form

$$\sum_{j,k=0}^{n-1} E(f^{(a)}(B_{j/n})f^{(b)}(B_{k/n}))\langle \varepsilon_{j/n}, \delta_{k/n}\rangle_{\mathfrak{H}} \prod_{i=1}^{3}\langle \varepsilon_{x_i/n}, \delta_{y_i/n}\rangle_{\mathfrak{H}}$$

(where $x_i$ and $y_i$ are for $j$ or $k$), and from Lemma 3.1(i), (iii), we have that

$$\begin{cases} \sup_{j,k=0,\ldots,n-1} \prod_{i=1}^{3} |\langle \varepsilon_{x_i/n}, \delta_{y_i/n} \rangle_{\mathfrak{H}}| \underset{n\to\infty}{=} O(n^{-3/2}), \\ \sum_{j,k=0}^{n-1} |\langle \varepsilon_{j/n}, \delta_{k/n} \rangle_{\mathfrak{H}}| \underset{n\to\infty}{=} O(n). \end{cases}$$

Still using the Malliavin integration by parts formula (2.5), we can bound $B_n$ as follows:

$$|B_n| \leq 8 \sum_{j,k=0}^{n-1} \sum_{a+b=2} |E(f^{(a)}(B_{j/n})f^{(b)}(B_{k/n}))$$

$$\times \langle \varepsilon_{j/n}^{\otimes a} \widetilde{\otimes} \varepsilon_{k/n}^{\otimes b}, \delta_{j/n} \otimes \delta_{k/n} \rangle_{\mathfrak{H}^{\otimes 2}} \langle \delta_{j/n}, \delta_{k/n} \rangle_{\mathfrak{H}}|$$

$$\leq Cn^{-1} \sum_{j,k=0}^{n-1} |\langle \delta_{j/n}, \delta_{k/n} \rangle_{\mathfrak{H}}| \qquad \text{by Lemma 3.1(i) and under (H}_4)$$

$$= Cn^{-3/2} \sum_{j,k=0}^{n-1} |\rho(j-k)|$$

$$\leq Cn^{-1/2} \sum_{r=-\infty}^{\infty} |\rho(r)| \underset{n\to\infty}{=} O(n^{-1/2}),$$

where

(3.4) $$\rho(r) := \sqrt{|r+1|} + \sqrt{|r-1|} - 2\sqrt{|r|}, \qquad r \in \mathbb{Z}.$$

Observe that the series $\sum_{r=-\infty}^{\infty} |\rho(r)|$ is convergent since $|\rho(r)| \underset{|r|\to\infty}{\sim} \frac{1}{2}|r|^{-3/2}$.

Finally, we consider the term $C_n$.

$$\begin{aligned} C_n &= \frac{1}{2n} \sum_{j,k=0}^{n-1} E(f(B_{j/n})f(B_{k/n}))\rho^2(j-k) \\ &= \frac{1}{2n} \sum_{r=-\infty}^{\infty} \sum_{j=0\vee -r}^{(n-1)\wedge(n-1-r)} E(f(B_{j/n})f(B_{(j+r)/n}))\rho^2(r) \\ &\underset{n\to\infty}{\longrightarrow} \frac{1}{2} \int_0^1 E(f^2(B_s))\, ds \sum_{r=-\infty}^{\infty} \rho^2(r) \\ &= C_{1/4}^2 \int_0^1 E(f^2(B_s))\, ds. \end{aligned}$$

The desired convergence (3.2) follows.

STEP 2. Since the sequence $(G_n)$ is bounded in $L^2$, the sequence $(G_n, (B_t)_{t \in [0,1]})$ is tight in $\mathbb{R} \times \mathscr{C}([0,1])$. Assume that $(G_\infty, (B_t)_{t \in [0,1]})$ denotes the limit in law of a certain subsequence of $(G_n, (B_t)_{t \in [0,1]})$, denoted again by $(G_n, (B_t)_{t \in [0,1]})$.

We have to prove that

$$G_\infty \stackrel{\text{Law}}{=} C_{1/4} \int_0^1 f(B_s) \, dW_s + \frac{1}{4} \int_0^1 f''(B_s) \, ds,$$

where $W$ denotes a standard Brownian motion independent of $B$, or, equivalently, that

(3.5)
$$E(e^{i\lambda G_\infty} | (B_t)_{t \in [0,1]}) = \exp\left\{ i \frac{\lambda}{4} \int_0^1 f''(B_s) \, ds - \frac{\lambda^2}{2} C_{1/4}^2 \int_0^1 f^2(B_s) \, ds \right\}.$$

This will be done by showing that for every random variable $\xi$ of the form (2.4) and every real number $\lambda$, we have

(3.6) $$\lim_{n \to \infty} \phi'_n(\lambda) = E\left\{ e^{i\lambda G_\infty} \xi \left( \frac{i}{4} \int_0^1 f''(B_s) \, ds - \lambda C_{1/4}^2 \int_0^1 f^2(B_s) \, ds \right) \right\},$$

where

$$\phi'_n(\lambda) := \frac{d}{d\lambda} E(e^{i\lambda G_n} \xi) = iE(G_n e^{i\lambda G_n} \xi), \qquad n \geq 1.$$

Let us make precise this argument. Because $(G_\infty, (B_t)_{t \in [0,1]})$ is the limit in law of $(G_n, (B_t)_{t \in [0,1]})$ and $(G_n)$ is bounded in $L^2$, we have that

$$E(G_\infty \xi e^{i\lambda G_\infty}) = \lim_{n \to \infty} E(G_n \xi e^{i\lambda G_n}) \qquad \forall \lambda \in \mathbb{R}$$

for every $\xi$ of the form (2.4). Furthermore, because convergence (3.6) holds for every $\xi$ of the form (2.4), the conditional characteristic function $\lambda \mapsto E(e^{i\lambda G_\infty} | (B_t)_{t \in [0,1]})$ satisfies the following linear ordinary differential equation:

$$\frac{d}{d\lambda} E(e^{i\lambda G_\infty} | (B_t)_{t \in [0,1]})$$
$$= E(e^{i\lambda G_\infty} | (B_t)_{t \in [0,1]}) \left[ \frac{i}{4} \int_0^1 f''(B_s) \, ds - \lambda C_{1/4}^2 \int_0^1 f^2(B_s) \, ds \right].$$

By solving it, we obtain (3.5), which yields the desired conclusion.

Thus it remains to show (3.6). By the duality between the derivative and divergence operators, we have

(3.7) $$E(f(B_{k/n}) I_2(\delta_{k/n}^{\otimes 2}) e^{i\lambda G_n} \xi) = E(\langle D^2(f(B_{k/n}) e^{i\lambda G_n} \xi), \delta_{k/n}^{\otimes 2} \rangle_{\mathfrak{H}^{\otimes 2}}).$$

The first and second derivatives of $f(B_{k/n})e^{i\lambda G_n}\xi$ are given by

$$D(f(B_{k/n})e^{i\lambda G_n}\xi) = f'(B_{k/n})e^{i\lambda G_n}\xi\varepsilon_{k/n} + i\lambda f(B_{k/n})e^{i\lambda G_n}\xi DG_n$$
$$+ f(B_{k/n})e^{i\lambda G_n}D\xi$$

and

$$D^2(f(B_{k/n})e^{i\lambda G_n}\xi)$$
$$= f''(B_{k/n})e^{i\lambda G_n}\xi\varepsilon_{k/n}^{\otimes 2} + 2i\lambda f'(B_{k/n})e^{i\lambda G_n}\xi(\varepsilon_{k/n}\widetilde{\otimes}DG_n)$$
$$+ 2f'(B_{k/n})e^{i\lambda G_n}(\varepsilon_{k/n}\widetilde{\otimes}D\xi) - \lambda^2 f(B_{k/n})e^{i\lambda G_n}\xi DG_n^{\otimes 2}$$
$$+ 2i\lambda f(B_{k/n})e^{i\lambda G_n}(DG_n\widetilde{\otimes}D\xi)$$
$$+ i\lambda f(B_{k/n})e^{i\lambda G_n}\xi D^2 G_n + f(B_{k/n})e^{i\lambda G_n}D^2\xi.$$

Hence allowing for expectation and multiplying by $\delta_{k/n}^{\otimes 2}$ yields

$$(3.8) \quad \begin{aligned} &E(\langle D^2(f(B_{k/n})e^{i\lambda G_n}\xi), \delta_{k/n}^{\otimes 2}\rangle_{\mathfrak{H}^{\otimes 2}}) \\ &= E(f''(B_{k/n})e^{i\lambda G_n}\xi)\langle\varepsilon_{k/n},\delta_{k/n}\rangle_{\mathfrak{H}}^2 \\ &\quad + 2i\lambda E(f'(B_{k/n})e^{i\lambda G_n}\xi\langle DG_n,\delta_{k/n}\rangle_{\mathfrak{H}})\langle\varepsilon_{k/n},\delta_{k/n}\rangle_{\mathfrak{H}} \\ &\quad + 2E(f'(B_{k/n})e^{i\lambda G_n}\langle D\xi,\delta_{k/n}\rangle_{\mathfrak{H}})\langle\varepsilon_{k/n},\delta_{k/n}\rangle_{\mathfrak{H}} \\ &\quad - \lambda^2 E(f(B_{k/n})e^{i\lambda G_n}\xi\langle DG_n,\delta_{k/n}\rangle_{\mathfrak{H}}^2) \\ &\quad + 2i\lambda E(f(B_{k/n})e^{i\lambda G_n}\langle D\xi,\delta_{k/n}\rangle_{\mathfrak{H}}\langle DG_n,\delta_{k/n}\rangle_{\mathfrak{H}}) \\ &\quad + i\lambda E(f(B_{k/n})e^{i\lambda G_n}\xi\langle D^2 G_n,\delta_{k/n}^{\otimes 2}\rangle_{\mathfrak{H}^{\otimes 2}}) \\ &\quad + E(f(B_{k/n})e^{i\lambda G_n}\langle D^2\xi,\delta_{k/n}^{\otimes 2}\rangle_{\mathfrak{H}^{\otimes 2}}). \end{aligned}$$

We also need explicit expressions for $\langle DG_n,\delta_{k/n}\rangle_{\mathfrak{H}}$ and for $\langle D^2 G_n,\delta_{k/n}^{\otimes 2}\rangle_{\mathfrak{H}^{\otimes 2}}$. By differentiating $G_n$, we obtain

$$(3.9) \quad DG_n = \sum_{l=0}^{n-1}[f'(B_{l/n})I_2(\delta_{l/n}^{\otimes 2})\varepsilon_{l/n} + 2f(B_{l/n})\Delta B_{l/n}\delta_{l/n}].$$

As a consequence,

$$(3.10) \quad \begin{aligned} \langle DG_n,\delta_{k/n}\rangle_{\mathfrak{H}} &= \sum_{l=0}^{n-1} f'(B_{l/n})I_2(\delta_{l/n}^{\otimes 2})\langle\varepsilon_{l/n},\delta_{k/n}\rangle_{\mathfrak{H}} \\ &\quad + 2\sum_{l=0}^{n-1} f(B_{l/n})\Delta B_{l/n}\langle\delta_{l/n},\delta_{k/n}\rangle_{\mathfrak{H}}. \end{aligned}$$

Also

$$D^2 G_n = \sum_{l=0}^{n-1} [f''(B_{l/n}) I_2(\delta_{l/n}^{\otimes 2}) \varepsilon_{l/n}^{\otimes 2}$$
$$+ 4f'(B_{l/n}) \Delta B_{l/n} (\varepsilon_{l/n} \widetilde{\otimes} \delta_{l/n}) + 2f(B_{l/n}) \delta_{l/n}^{\otimes 2}],$$

and, as a consequence,

$$\langle D^2 G_n, \delta_{k/n}^{\otimes 2} \rangle_{\mathfrak{H}^{\otimes 2}} = \sum_{l=0}^{n-1} [f''(B_{l/n}) I_2(\delta_{l/n}^{\otimes 2}) \langle \varepsilon_{l/n}, \delta_{k/n} \rangle_{\mathfrak{H}}^2$$
(3.11)
$$+ 4f'(B_{l/n}) \Delta B_{l/n} \langle \varepsilon_{l/n}, \delta_{k/n} \rangle_{\mathfrak{H}} \langle \delta_{l/n}, \delta_{k/n} \rangle_{\mathfrak{H}}$$
$$+ 2f(B_{l/n}) \langle \delta_{l/n}, \delta_{k/n} \rangle_{\mathfrak{H}}^2].$$

Substituting (3.11) into (3.8) yields the following decomposition for $\phi'_n(\lambda) = iE(G_n e^{i\lambda G_n} \xi)$:

$$\phi'_n(\lambda) = -2\lambda \sum_{k,l=0}^{n-1} E(f(B_{k/n}) f(B_{l/n}) e^{i\lambda G_n} \xi) \langle \delta_{l/n}, \delta_{k/n} \rangle_{\mathfrak{H}}^2$$
(3.12)
$$+ i \sum_{k=0}^{n-1} E(f''(B_{k/n}) e^{i\lambda G_n} \xi) \langle \varepsilon_{k/n}, \delta_{k/n} \rangle_{\mathfrak{H}}^2 + i \sum_{k=0}^{n-1} r_{k,n},$$

where $r_{k,n}$ is given by

$$r_{k,n} = 2i\lambda E(f'(B_{k/n}) e^{i\lambda G_n} \xi \langle DG_n, \delta_{k/n} \rangle_{\mathfrak{H}}) \langle \varepsilon_{k/n}, \delta_{k/n} \rangle_{\mathfrak{H}}$$
$$+ 2E(f'(B_{k/n}) e^{i\lambda G_n} \langle D\xi, \delta_{k/n} \rangle_{\mathfrak{H}}) \langle \varepsilon_{k/n}, \delta_{k/n} \rangle_{\mathfrak{H}}$$
$$- \lambda^2 E(f(B_{k/n}) e^{i\lambda G_n} \xi \langle DG_n, \delta_{k/n} \rangle_{\mathfrak{H}}^2)$$
(3.13)
$$+ 2i\lambda E(f(B_{k/n}) e^{i\lambda G_n} \langle D\xi, \delta_{k/n} \rangle_{\mathfrak{H}} \langle DG_n, \delta_{k/n} \rangle_{\mathfrak{H}})$$
$$+ i\lambda \sum_{l=0}^{n-1} E(f(B_{k/n}) e^{i\lambda G_n} \xi f''(B_{l/n}) I_2(\delta_{l/n}^{\otimes 2})) \langle \varepsilon_{l/n}, \delta_{k/n} \rangle_{\mathfrak{H}}^2$$
$$+ 4i\lambda \sum_{l=0}^{n-1} E(f(B_{k/n}) e^{i\lambda G_n} \xi f'(B_{l/n}) \Delta B_{l/n}) \langle \varepsilon_{l/n}, \delta_{k/n} \rangle_{\mathfrak{H}} \langle \delta_{l/n}, \delta_{k/n} \rangle_{\mathfrak{H}}$$
$$+ E(f(B_{k/n}) e^{i\lambda G_n} \langle D^2 \xi, \delta_{k/n}^{\otimes 2} \rangle_{\mathfrak{H}^{\otimes 2}}) = \sum_{j=1}^{7} R_{k,n}^{(j)}.$$

Remark that the first sum in the right-hand side of (3.12) is very similar to the quantity $C_n$ presented in Step 1. In fact, similar computations give

$$\lim_{n\to\infty} -2\lambda \sum_{k,l=0}^{n-1} \mathbb{E}[f(B_{k/n})f(B_{l/n})e^{i\lambda G_n}\xi]\langle \delta_{l/n}, \delta_{k/n}\rangle_{\mathfrak{H}}^2$$
(3.14)
$$= -C_{1/4}^2 \lambda \int_0^1 E(f^2(B_s)e^{i\lambda G_\infty}\xi)\,ds.$$

Furthermore, the second term of (3.12) is very similar to $E(G_n)$. In fact, using the arguments presented in Step 1, we obtain here that

$$\lim_{n\to\infty} i \sum_{k=0}^{n-1} E(f''(B_{k/n})e^{i\lambda G_n}\xi)\langle \varepsilon_{k/n}, \delta_{k/n}\rangle_{\mathfrak{H}}^2$$
(3.15)
$$= \frac{i}{4} \int_0^1 E(f''(B_s)e^{i\lambda G_\infty}\xi)\,ds.$$

Consequently, (3.6) will be shown if we prove that $\lim_{n\to\infty}\sum_{k=0}^{n-1} r_{k,n} = 0$. This will be done in several steps.

STEP 3. In this step, we state and prove some estimates which are crucial in the rest of the proof. First, we will show that

(3.16) $\quad |E(f'(B_{k/n})f'(B_{l/n})e^{i\lambda G_n}\xi I_2(\delta_{l/n}^{\otimes 2}))| \leq \dfrac{C}{n} \quad$ for any $0 \leq k, l \leq n-1$.

Then we will prove that

$$|E(f(B_{k/n})f'(B_{j/n})f'(B_{l/n})e^{i\lambda G_n}\xi I_4(\delta_{j/n}^{\otimes 2} \otimes \delta_{l/n}^{\otimes 2}))| \leq \frac{C}{n^2}$$
(3.17)
$$\text{for any } 0 \leq k, j, l \leq n-1.$$

PROOF OF (3.16). Let $\zeta_{\xi,k,n}$ denote any random variable of the form $f^{(a)}(B_{k/n})f^{(b)}(B_{l/n})e^{i\lambda G_n}\xi$ with $a$ and $b$ two positive integers less or equal to four. From the Malliavin integration by parts formula (2.5) we have

$$E(f'(B_{k/n})f'(B_{l/n})e^{i\lambda G_n}\xi I_2(\delta_{l/n}^{\otimes 2}))$$
$$= E(\langle D^2(f'(B_{k/n})f'(B_{l/n})e^{i\lambda G_n}\xi), \delta_{l/n}^{\otimes 2}\rangle_{\mathfrak{H}^{\otimes 2}}).$$

When computing the right-hand side, three types of terms appear. First, we have some terms of the form,

(3.18) $\quad \begin{cases} E(\zeta_{\xi,k,n})\langle \varepsilon_{k/n}, \delta_{l/n}\rangle_{\mathfrak{H}}^2, \text{ or} \\ E(\zeta_{\xi,k,n}\langle D\xi, \delta_{l/n}\rangle_{\mathfrak{H}})\langle \varepsilon_{k/n}, \delta_{l/n}\rangle_{\mathfrak{H}}, \text{ or} \\ E(\zeta_{\xi,k,n}\langle D^2\xi, \delta_{l/n}^{\otimes 2}\rangle_{\mathfrak{H}^{\otimes 2}}), \end{cases}$

where $D\xi$ and $D^2\xi$ are given by,

$$\begin{cases} D\xi = \sum_{i=1}^{m} \dfrac{\partial \psi}{\partial x_i}(B_{t_1},\ldots,B_{t_m})\varepsilon_{t_i}, \\ D^2\xi = \sum_{i,j=1}^{m} \dfrac{\partial^2 \psi}{\partial x_j \, \partial x_i}(B_{t_1},\ldots,B_{t_m})\varepsilon_{t_j} \otimes \varepsilon_{t_i}. \end{cases}$$

From Lemma 3.1(i) and under (H$_4$), we have that each of the three terms in (3.18) is less or equal to $Cn^{-1}$. The second type of term we have to deal with is

(3.19) $$\begin{cases} E(\zeta_{\xi,k,n}\langle DG_n, \delta_{l/n}\rangle_{\mathfrak{H}})\langle \varepsilon_{k/n}, \delta_{l/n}\rangle_{\mathfrak{H}}, \text{ or} \\ E(\zeta_{\xi,k,n}\langle DG_n, \delta_{l/n}\rangle_{\mathfrak{H}}\langle D\xi, \delta_{l/n}\rangle_{\mathfrak{H}}). \end{cases}$$

By the Cauchy–Schwarz inequality, under (H$_4$) and by using (4.20) in the third version of [16], that is,

$$E(\langle DG_n, \delta_{l/n}\rangle_{\mathfrak{H}}^2) \leq Cn^{-1},$$

we have that both expressions in (3.19) are also less or equal to $Cn^{-1}$.

The last type of term which has to be taken into account is the term

$$-\lambda^2 E(f'(B_{k/n})f'(B_{l/n})e^{i\lambda G_n}\xi\langle D^2 G_n, \delta_{k/n}^{\otimes 2}\rangle_{\mathfrak{H}^{\otimes 2}}).$$

Again, by using the Cauchy–Schwarz inequality and the estimate

$$E(\langle D^2 G_n, \delta_{k/n}^{\otimes 2}\rangle_{\mathfrak{H}^{\otimes 2}}^2) \leq Cn^{-2}$$

(which can be obtained by mimicking the proof of (4.20) in the third version of [16]), we can conclude that

$$|-\lambda^2 E(f'(B_{k/n})f'(B_{l/n})e^{i\lambda G_n}\xi\langle D^2 G_n, \delta_{k/n}^{\otimes 2}\rangle_{\mathfrak{H}^{\otimes 2}})| \leq \frac{C}{n}.$$

As a consequence (3.16) is shown.

PROOF OF (3.17). By the Malliavin integration by parts formula (2.5), we have

$$E(\zeta_{\xi,k,n}f'(B_{j/n})f'(B_{l/n})I_4(\delta_{j/n}^{\otimes 2} \otimes \delta_{l/n}^{\otimes 2}))$$
$$= E(\langle D^4(\zeta_{\xi,k,n}f'(B_{j/n})f'(B_{l/n})), \delta_{j/n}^{\otimes 2} \otimes \delta_{l/n}^{\otimes 2}\rangle_{\mathfrak{H}^{\otimes 4}}).$$

When computing the right-hand side, we have to deal with the same type of term as in the proof of (3.16), plus two additional types of terms containing

$$E(\langle D^3 G_n, \delta_{j/n}^{\otimes 2} \otimes \delta_{l/n}\rangle_{\mathfrak{H}^{\otimes 3}}^2) \quad \text{and} \quad E(\langle D^4 G_n, \delta_{j/n}^{\otimes 2} \otimes \delta_{l/n}^{\otimes 2}\rangle_{\mathfrak{H}^{\otimes 4}}^2).$$

In fact, by mimicking the proof of (4.20) in the third version of [16], we can obtain the following bounds:

$$E(\langle D^3 G_n, \delta_{j/n}^{\otimes 2} \otimes \delta_{l/n}\rangle_{\mathfrak{H}^{\otimes 3}}^2) \leq Cn^{-3} \quad \text{and} \quad E(\langle D^4 G_n, \delta_{j/n}^{\otimes 2} \otimes \delta_{l/n}^{\otimes 2}\rangle_{\mathfrak{H}^{\otimes 4}}^2) \leq Cn^{-4}.$$

This allows us to obtain (3.17).

STEP 4. We compute the terms corresponding to $R_{k,n}^{(1)}$, $R_{k,n}^{(4)}$ and $R_{k,n}^{(6)}$ in (3.13). The derivative $DG_n$ is given by (3.9) so that

$$\sum_{k=0}^{n-1} R_{k,n}^{(1)} = 2i\lambda \sum_{k,l=0}^{n-1} E(f'(B_{k/n})f'(B_{l/n})e^{i\lambda G_n}\xi I_2(\delta_{l/n}^{\otimes 2}))\langle \varepsilon_{l/n}, \delta_{k/n}\rangle_{\mathfrak{H}} \langle \varepsilon_{k/n}, \delta_{k/n}\rangle_{\mathfrak{H}}$$

$$+ 2\sum_{k,l=0}^{n-1} E(f'(B_{k/n})f(B_{l/n})e^{i\lambda G_n}\xi \Delta B_{l/n})\langle \delta_{l/n}, \delta_{k/n}\rangle_{\mathfrak{H}} \langle \varepsilon_{k/n}, \delta_{k/n}\rangle_{\mathfrak{H}}$$

$$= T_1^{(1)} + T_2^{(1)}.$$

From (3.16), Lemma 3.1(i), (iii) and under $(H_4)$, we have that

$$|T_1^{(1)}| \leq Cn^{-3/2} \sum_{k,l=0}^{n-1} |\langle \varepsilon_{l/n}, \delta_{k/n}\rangle_{\mathfrak{H}}| \leq Cn^{-1/2}.$$

For $T_2^{(1)}$, remark first that the Cauchy–Schwarz inequality and hypothesis $(H_4)$ yields

(3.20) $\qquad |E(f'(B_{k/n})e^{i\lambda G_n}\xi f(B_{l/n})\Delta B_{l/n})| \leq Cn^{-1/4}.$

Thus by Lemma 3.1(i),

$$|T_2^{(1)}| \leq Cn^{-3/4} \sum_{k,l=0}^{n-1} |\langle \delta_{l/n}, \delta_{k/n}\rangle_{\mathfrak{H}}| = Cn^{-5/4} \sum_{k,l=0}^{n-1} |\rho(k-l)|$$

$$\leq Cn^{-1/4} \sum_{r=-\infty}^{\infty} |\rho(r)| = Cn^{-1/4},$$

where $\rho$ has been defined in (3.4).

The term corresponding to $R_{k,n}^{(4)}$ is very similar to $R_{k,n}^{(1)}$. Indeed, by (3.9), we have

$$\sum_{k=0}^{n-1} R_{k,n}^{(4)} = 2i\lambda \sum_{i=1}^{m} \sum_{k,l=0}^{n-1} E\left(f(B_{k/n})f'(B_{l/n})e^{i\lambda G_n} \frac{\partial \psi}{\partial x_i}(B_{t_1},\ldots,B_{t_m})I_2(\delta_{l/n}^{\otimes 2})\right)$$

$$\times \langle \varepsilon_{l/n}, \delta_{k/n}\rangle_{\mathfrak{H}} \langle \varepsilon_{t_i}, \delta_{k/n}\rangle_{\mathfrak{H}}$$

$$+ 4i\lambda \sum_{i=1}^{m} \sum_{k,l=0}^{n-1} E\left(f(B_{k/n})f(B_{l/n})e^{i\lambda G_n}\Delta B_{l/n} \frac{\partial \psi}{\partial x_i}(B_{t_1},\ldots,B_{t_m})\right)$$

$$\times \langle \delta_{l/n}, \delta_{k/n}\rangle_{\mathfrak{H}} \langle \varepsilon_{t_i}, \delta_{k/n}\rangle_{\mathfrak{H}}$$

$$= T_1^{(4)} + T_2^{(4)}$$

and we can proceed for $T_i^{(4)}$ as for $T_i^{(1)}$.

The term corresponding to $R_{k,n}^{(6)}$ is very similar to $T_2^{(1)}$. More precisely, we have

$$\left|\sum_{k=0}^{n-1} R_{k,n}^{(6)}\right| \leq Cn^{-3/4} \sum_{k,l=0}^{n-1} |\langle \delta_{l/n}, \delta_{k/n}\rangle_{\mathfrak{H}}| = Cn^{-5/4} \sum_{k,l=0}^{n-1} |\rho(k-l)|$$

$$\leq Cn^{-1/4} \sum_{r=-\infty}^{\infty} |\rho(r)| = Cn^{-1/4}.$$

STEP 5. Estimation of $R_{k,n}^{(3)}$. Let $\zeta_{\xi,k,n} := \lambda^2 f(B_{k/n})e^{i\lambda G_n}\xi$. Using (3.9), we have

$$\langle DG_n, \delta_{k/n}\rangle_{\mathfrak{H}}^2 = \sum_{j,l=0}^{n-1} f'(B_{l/n})f'(B_{j/n}) I_2(\delta_{l/n}^{\otimes 2}) I_2(\delta_{j/n}^{\otimes 2}) \langle \varepsilon_{j/n}, \delta_{k/n}\rangle_{\mathfrak{H}} \langle \varepsilon_{l/n}, \delta_{k/n}\rangle_{\mathfrak{H}}$$

$$+ \sum_{j,l=0}^{n-1} f(B_{j/n}) f(B_{l/n}) \Delta B_{j/n} \Delta B_{l/n} \langle \delta_{j/n}, \delta_{k/n}\rangle_{\mathfrak{H}} \langle \delta_{l/n}, \delta_{k/n}\rangle_{\mathfrak{H}},$$

and, consequently,

$$\left|\sum_{k=0}^{n-1} R_{k,n}^3\right| \leq \sum_{k=0}^{n-1} |E(\zeta_{\xi,k,n}\langle DG_n, \delta_{k/n}\rangle_{\mathfrak{H}}^2)|$$

$$\leq 2 \sum_{k,j,l=0}^{n-1} |E(\zeta_{\xi,k,n} f'(B_{j/n}) f'(B_{l/n}) I_2(\delta_{j/n}^{\otimes 2}) I_2(\delta_{l/n}^{\otimes 2}))$$

$$\times \langle \varepsilon_{j/n}, \delta_{k/n}\rangle_{\mathfrak{H}} \langle \varepsilon_{l/n}, \delta_{k/n}\rangle_{\mathfrak{H}}|$$

$$+ 8 \sum_{k,j,l=0}^{n-1} |E(\zeta_{\xi,k,n} f(B_{j/n}) f(B_{l/n}) \Delta B_{j/n} \Delta B_{l/n})$$

$$\times \langle \delta_{j/n}, \delta_{k/n}\rangle_{\mathfrak{H}} \langle \delta_{l/n}, \delta_{k/n}\rangle_{\mathfrak{H}}|.$$

Using the product formula (3.3), we have

$$\left|\sum_{k=0}^{n-1} R_{k,n}^3\right| \leq 2 \sum_{k,j,l=0}^{n-1} |E(\zeta_{\xi,k,n} f'(B_{j/n}) f'(B_{l/n}) I_4(\delta_{j/n}^{\otimes 2} \otimes \delta_{l/n}^{\otimes 2}))|$$

$$\times |\langle \varepsilon_{j/n}, \delta_{k/n}\rangle_{\mathfrak{H}} \langle \varepsilon_{l/n}, \delta_{k/n}\rangle_{\mathfrak{H}}|$$

$$+ 8 \sum_{k,j,l=0}^{n-1} |E(\zeta_{\xi,k,n} f'(B_{j/n}) f'(B_{l/n}) I_2(\delta_{j/n} \otimes \delta_{l/n}))|$$

$$\times |\langle \delta_{j/n}, \delta_{l/n}\rangle_{\mathfrak{H}} \langle \varepsilon_{j/n}, \delta_{k/n}\rangle_{\mathfrak{H}} \langle \varepsilon_{l/n}, \delta_{k/n}\rangle_{\mathfrak{H}}|$$

$$+ 4 \sum_{k,j,l=0}^{n-1} |E(\zeta_{\xi,k,n} f'(B_{j/n}) f'(B_{l/n}))|$$

$$\times |\langle \delta_{j/n}, \delta_{l/n} \rangle_{\mathfrak{H}}^2 \langle \varepsilon_{j/n}, \delta_{k/n} \rangle_{\mathfrak{H}} \langle \varepsilon_{l/n}, \delta_{k/n} \rangle_{\mathfrak{H}}|$$

$$+ 8 \sum_{k,j,l=0}^{n-1} |E(\zeta_{\xi,k,n} f(B_{j/n}) f(B_{l/n}) \Delta B_{j/n} \Delta B_{l/n})|$$

$$\times |\langle \delta_{j/n}, \delta_{k/n} \rangle_{\mathfrak{H}} \langle \delta_{l/n}, \delta_{k/n} \rangle_{\mathfrak{H}}|$$

$$= \sum_{i=1}^{4} T_i^{(3)}.$$

From (3.17), we have

$$|T_1^{(3)}| \leq C n^{-1/2} \sum_{k,j,l=0}^{n-1} |E(\zeta_{\xi,k,n} f'(B_{j/n}) f'(B_{l/n}) I_4(\delta_{j/n}^{\otimes 2} \otimes \delta_{l/n}^{\otimes 2}))| |\langle \varepsilon_{j/n}, \delta_{k/n} \rangle_{\mathfrak{H}}|$$

$$\leq C n^{-5/2} n^2 \sup_{j=0,\ldots,n-1} \sum_{k=0}^{n-1} |\langle \varepsilon_{j/n}, \delta_{k/n} \rangle_{\mathfrak{H}}| \leq C n^{-1/2} \quad \text{by Lemma 3.1(ii).}$$

Now let us consider $T_2^{(3)}$. Using (3.16) and Lemma 3.1(ii), we deduce that

$$|T_2^{(3)}| \leq C n^{-3/2} \sum_{j,l=0}^{n-1} |\langle \delta_{j/n}, \delta_{l/n} \rangle_{\mathfrak{H}}| \sup_{j=0,\ldots,n-1} \sum_{k=0}^{n-1} |\langle \varepsilon_{j/n}, \delta_{k/n} \rangle_{\mathfrak{H}}|$$

$$\leq C n^{-1/2} \sum_{r=-\infty}^{\infty} |\rho(r)| = C n^{-1/2}.$$

For $T_3^{(3)}$, we have

$$|T_3^{(3)}| \leq C n^{-1/2} \sum_{j,l=0}^{n-1} \langle \delta_{j/n}, \delta_{l/n} \rangle_{\mathfrak{H}}^2 \sup_{j=0,\ldots,n-1} \sum_{k=0}^{n-1} |\langle \varepsilon_{j/n}, \delta_{k/n} \rangle_{\mathfrak{H}}|$$

$$\leq C n^{-1/2} \sum_{r=-\infty}^{\infty} \rho^2(r) = C n^{-1/2}.$$

Finally, by the Cauchy–Schwarz inequality and under (H$_4$), we have

$$|E(\zeta_{\xi,k,n} f(B_{j/n}) f(B_{l/n}) \Delta B_{j/n} \Delta B_{l/n})| \leq C n^{-1/2}.$$

Consequently,

$$|T_4^{(3)}| \leq C n^{-1/2} \sum_{k,j,l=0}^{n-1} |\langle \delta_{j/n}, \delta_{k/n} \rangle_{\mathfrak{H}} \langle \delta_{k/n}, \delta_{l/n} \rangle_{\mathfrak{H}}|$$

$$\leq Cn^{-3/2} \sum_{k,j,l=0}^{n-1} |\rho(k-l)\rho(k-j)|$$

$$\leq Cn^{-1/2} \left( \sum_{r=-\infty}^{\infty} |\rho(r)| \right)^2$$

$$= Cn^{-1/2}.$$

STEP 6. Estimation of $R_{k,n}^{(5)}$. From (3.16) and Lemma 3.1(iii), we have

$$\left| \sum_{k=0}^{n-1} R_{k,n}^{(5)} \right| \leq Cn^{-3/2} \sum_{k,l=0}^{n-1} |\langle \varepsilon_{l/n}, \delta_{k/n} \rangle_{\mathfrak{H}}| \leq Cn^{-1/2}.$$

STEP 7. Estimation of $R_{k,n}^{(2)}$ and $R_{k,n}^{(7)}$. We recall that

$$0 \leq \sqrt{x+1} - \sqrt{x} \leq 1 \qquad \text{for any } x \geq 0.$$

Thus under (H$_4$) and using Lemma 3.1, we have,

$$\left| \sum_{k=0}^{n-1} R_{k,n}^{(2)} \right| \leq 2 \sum_{i=1}^{m} \sum_{k=0}^{n-1} \left| E\left( f'(B_{k/n}) e^{i\lambda G_n} \frac{\partial \psi}{\partial x_i}(B_{t_1}, \ldots, B_{t_m}) \right) \right.$$
$$\left. \times \langle \varepsilon_{t_i}, \delta_{k/n} \rangle_{\mathfrak{H}} \langle \varepsilon_{k/n}, \delta_{k/n} \rangle_{\mathfrak{H}} \right|$$
$$\leq C(f, \psi) n^{-1/2} \sup_{t \in [0,1]} \sum_{k=0}^{n-1} |\langle \varepsilon_t, \delta_{k/n} \rangle_{\mathfrak{H}}| \leq Cn^{-1/2}.$$

Similarly, the following bound holds:

$$\left| \sum_{k=0}^{n-1} R_{k,n}^{(7)} \right| \leq \sum_{i,j=1}^{m} \sum_{k=0}^{n-1} \left| E\left( f(B_{k/n}) e^{i\lambda G_n} \frac{\partial^2 \psi}{\partial x_j \partial x_i}(B_{t_1}, \ldots, B_{t_m}) \right) \right.$$
$$\left. \times \langle \varepsilon_{t_i}, \delta_{k/n} \rangle_{\mathfrak{H}} \langle \varepsilon_{t_j}, \delta_{k/n} \rangle_{\mathfrak{H}} \right|$$
$$\leq Cn^{-1/2}.$$

The proof of Theorem 3.2 is complete. □

**4. Proof of Theorem 1.2.** Let $B = B^{1/4}$ be a fractional Brownian motion with the Hurst index $H = 1/4$. Moreover, we continue to note $\Delta B_{k/n}$ (resp. $\delta_{k/n}$; $\varepsilon_{k/n}$) instead of $B_{(k+1)/n} - B_{k/n}$ (resp. $\mathbf{1}_{[k/n,(k+1)/n]}$; $\mathbf{1}_{[0,k/n]}$). The aim of this section is to prove Theorem 1.2, or equivalently,

THEOREM 4.1 (Itô's formula). *Let $f:\mathbb{R} \to \mathbb{R}$ be a function verifying* (H$_9$). *Then*

$$\int_0^1 f'(B_s) \, d^\star B_s := \lim_{n \to \infty} \sum_{k=1}^{\lfloor n/2 \rfloor} f'(B_{(2k-1)/n})(B_{(2k)/n} - B_{(2k-2)/n})$$

*exists in law.*

*Moreover, we have*

$$\int_0^1 f'(B_s) \, d^\star B_s \stackrel{\text{Law}}{=} f(B_1) - f(0) - \frac{\kappa}{2} \int_0^1 f''(B_s) \, dW_s$$

*with $\kappa$ defined by*

$$(4.1) \qquad \kappa = \sqrt{2 + \sum_{r=1}^\infty (-1)^r \rho^2(r)} \approx 1290$$

*[recall the definition (3.4) of $\rho$] and where $W$ denotes a standard Brownian motion independent of $B$.*

PROOF. In [23] [identity (1.6)], it is proved that

$$\sum_{k=1}^{\lfloor n/2 \rfloor} f'(B_{(2k-1)/n})(B_{(2k)/n} - B_{(2k-2)/n})$$

$$\approx f(B_1) - f(0)$$

$$- \frac{1}{2} \sum_{k=1}^{\lfloor n/2 \rfloor} f''(B_{(2k-1)/n})[(\Delta B_{(2k-1)/n})^2 - (\Delta B_{(2k-2)/n})^2]$$

$$- \frac{1}{6} \sum_{j=1}^{\lfloor n/2 \rfloor} f'''(B_{(2j-1)/n})[(\Delta B_{(2j-2)/n})^3 + (\Delta B_{(2j-1)/n})^3],$$

where "$\approx$" means the difference goes to zero in $L^2$. Therefore, Theorem 4.1 is a direct consequence of Lemmas 4.2 and 4.3 below. $\square$

LEMMA 4.2. *Let $f:\mathbb{R} \to \mathbb{R}$ be a function verifying* (H$_6$). *Then*

$$(4.2) \qquad \sum_{j=1}^{\lfloor n/2 \rfloor} f(B_{(2j-1)/n})[(\Delta B_{(2j-2)/n})^3 + (\Delta B_{(2j-1)/n})^3] \xrightarrow[n \to \infty]{L^2} 0.$$

PROOF. Let $H_3(x) = x^3 - 3x$ be the third Hermite polynomial. Using the relation between Hermite polynomials and multiple integrals (see Section 2),

remark that

$$(\Delta B_{(2j-2)/n})^3 + (\Delta B_{(2j-1)/n})^3$$
$$= n^{-3/4}\bigg[H_3(n^{1/4}\Delta B_{(2j-2)/n}) + H_3(n^{1/4}\Delta B_{(2j-1)/n})$$
$$+ \frac{3}{\sqrt{n}}(B_{(2j-2)/n} - B_{(2j)/n})\bigg]$$
$$= I_3(\delta^{\otimes 3}_{(2j-2)/n}) + I_3(\delta^{\otimes 3}_{(2j-1)/n}) + \frac{3}{\sqrt{n}}I_1(\mathbf{1}_{[(2j-2)/n,(2j)/n]}),$$

so that (4.2) can be shown by successively proving that

$$(4.3) \quad E\left|\frac{1}{\sqrt{n}}\sum_{j=1}^{\lfloor n/2 \rfloor} f(B_{(2j-1)/n})I_1(\mathbf{1}_{[(2j-2)/n,(2j)/n]})\right|^2 \xrightarrow[n\to+\infty]{} 0;$$

$$(4.4) \quad E\left|\sum_{j=1}^{\lfloor n/2 \rfloor} f(B_{(2j-1)/n})I_3(\delta^{\otimes 3}_{(2j-2)/n})\right|^2 \xrightarrow[n\to+\infty]{} 0;$$

$$(4.5) \quad E\left|\sum_{j=1}^{\lfloor n/2 \rfloor} f(B_{(2j-1)/n})I_3(\delta^{\otimes 3}_{(2j-1)/n})\right|^2 \xrightarrow[n\to+\infty]{} 0.$$

Let us first proceed with the proof of (4.3). We can write, using, in particular, (2.6),

$$E\left|\frac{1}{\sqrt{n}}\sum_{j=1}^{\lfloor n/2 \rfloor} f(B_{(2j-1)/n})I_1(\mathbf{1}_{[(2j-2)/n,(2j)/n]})\right|^2$$
$$= \frac{1}{n}\left|\sum_{j,k=1}^{\lfloor n/2 \rfloor} E\{f(B_{(2j-1)/n})f(B_{(2k-1)/n})\right.$$
$$\left.\times I_1(\mathbf{1}_{[(2j-2)/n,(2j)/n]})I_1(\mathbf{1}_{[(2k-2)/n,(2k)/n]})\}\right|$$
$$\leq \frac{1}{n}\sum_{j,k=1}^{\lfloor n/2 \rfloor}|E\{f(B_{(2j-1)/n})f(B_{(2k-1)/n})$$
$$\times I_2(\mathbf{1}_{[(2j-2)/n,(2j)/n]} \otimes \mathbf{1}_{[(2k-2)/n,(2k)/n]})\}|$$
$$+ \frac{1}{n\sqrt{n}}\sum_{j,k=1}^{\lfloor n/2 \rfloor}|E\{f(B_{(2j-1)/n})f(B_{(2k-1)/n})\rho(2j-2k)\}|$$

$$= \frac{2}{n} \sum_{a+b=2} \sum_{j,k=1}^{\lfloor n/2 \rfloor} |E\{f^{(a)}(B_{(2j-1)/n})f^{(b)}(B_{(2k-1)/n})\}|$$

$$\times |\langle \varepsilon_{(2j-1)/n}^{\otimes a} \otimes \varepsilon_{(2k-1)/n}^{\otimes b},$$

$$\mathbf{1}_{[(2j-2)/n,(2j)/n]} \otimes \mathbf{1}_{[(2k-2)/n,(2k)/n]} \rangle_{\mathfrak{H}^{\otimes 2}}|$$

$$+ \frac{1}{n\sqrt{n}} \sum_{j,k=1}^{\lfloor n/2 \rfloor} |E\{f(B_{(2j-1)/n})f(B_{(2k-1)/n})\rho(2j-2k)\}|.$$

But by Lemma 3.1(i), we have

$$|\langle \varepsilon_{(2j-1)/n}^{\otimes a} \otimes \varepsilon_{(2k-1)/n}^{\otimes b}, \mathbf{1}_{[(2j-2)/n,(2j)/n]} \otimes \mathbf{1}_{[(2k-2)/n,(2k)/n]} \rangle_{\mathfrak{H}^{\otimes 2}}|$$

$$\leq \frac{1}{\sqrt{n}} (|\langle \varepsilon_{(2j-1)/n}, \mathbf{1}_{[(2j-2)/n,(2j)/n]} \rangle_{\mathfrak{H}}| + |\langle \varepsilon_{(2k-1)/n}, \mathbf{1}_{[(2k-2)/n,(2k)/n]} \rangle_{\mathfrak{H}}|)$$

$$= \frac{1}{n} (\sqrt{2j} - \sqrt{2j-2} + \sqrt{2k} - \sqrt{2k-2}).$$

Thus under (H$_6$),

$$\sum_{a+b=2} \sum_{j,k=1}^{\lfloor n/2 \rfloor} |E\{f^{(a)}(B_{(2j-1)/n})f^{(b)}(B_{(2k-1)/n})\}|$$

$$\times |\langle \varepsilon_{(2j-1)/n}^{\otimes a} \otimes \varepsilon_{(2k-1)/n}^{\otimes b}, \mathbf{1}_{[(2j-2)/n,(2j)/n]} \otimes \mathbf{1}_{[(2k-2)/n,(2k)/n]} \rangle_{\mathfrak{H}^{\otimes 2}}|$$

$$= O(\sqrt{n}).$$

Moreover,

$$\sum_{j,k=1}^{\lfloor n/2 \rfloor} |E\{f(B_{(2j-1)/n})f(B_{(2k-1)/n})\rho(2j-2k)\}|$$

$$\leq C \sum_{j,k=1}^{\lfloor n/2 \rfloor} |\rho(2j-2k)| = O(n).$$

Finally, convergence (4.3) holds.

Now let us only proceed with the proof of (4.4), the proof of (4.5) being similar. We have

$$E \left| \sum_{j=1}^{\lfloor n/2 \rfloor} f(B_{(2j-1)/n}) I_3(\delta_{(2j-2)/n}^{\otimes 3}) \right|^2$$

$$= \sum_{j,k=1}^{\lfloor n/2 \rfloor} E\{f(B_{(2j-1)/n})f(B_{(2k-1)/n}) I_3(\delta_{(2j-2)/n}^{\otimes 3}) I_3(\delta_{(2k-2)/n}^{\otimes 3})\}$$

$$= \sum_{r=0}^{3} r! \binom{3}{r}^2 n^{-(3-r)/2} \sum_{j,k=1}^{\lfloor n/2 \rfloor} E\{f(B_{(2j-1)/n})f(B_{(2k-1)/n})$$

$$\times I_{2r}(\delta_{(2j-2)/n}^{\otimes r} \otimes \delta_{(2k-2)/n}^{\otimes r})\}$$

$$\times \rho^{3-r}(2j-2k).$$

To obtain (4.4), it is then sufficient to prove that, for every fixed $r \in \{0,1,2,3\}$, the quantities

$$R_n^{(r)} = n^{-(3-r)/2} \sum_{j,k=1}^{\lfloor n/2 \rfloor} E\{f(B_{(2j-1)/n})f(B_{(2k-1)/n})$$

$$\times I_{2r}(\delta_{(2j-2)/n}^{\otimes r} \otimes \delta_{(2k-2)/n}^{\otimes r})\}$$

$$\times \rho^{3-r}(2j-2k)$$

tend to zero as $n \to \infty$. We have, by Lemma 3.1(i) and under (H$_6$),

$$\sup_{j,k=1,\ldots,\lfloor n/2 \rfloor} |E\{f(B_{(2j-1)/n})f(B_{(2k-1)/n})I_{2r}(\delta_{(2j-2)/n}^{\otimes r} \otimes \delta_{(2k-2)/n}^{\otimes r})\}|$$

$$= \sup_{j,k=1,\ldots,\lfloor n/2 \rfloor} (2r)! \Big| \sum_{a+b=2r} E\{f^{(a)}(B_{(2j-1)/n})f^{(b)}(B_{(2k-1)/n})\}$$

$$\times \langle \varepsilon_{(2j-1)/n}^{\otimes a} \widetilde{\otimes} \varepsilon_{(2j-1)/n}^{\otimes b},$$

$$\mathbf{1}_{[(2j-2)/n,(2j/n)]}^{\otimes r} \otimes \mathbf{1}_{[(2k-2)/n,(2k/n)]}^{\otimes r} \rangle_{\mathfrak{H}^{\otimes 2}} \Big|$$

$$\leq C \sup_{j,k=1,\ldots,\lfloor n/2 \rfloor} \sup_{a+b=2r} |\langle \varepsilon_{(2j-1)/n}^{\otimes a} \widetilde{\otimes} \varepsilon_{(2j-1)/n}^{\otimes b},$$

$$\mathbf{1}_{[(2j-2)/n,(2j/n)]}^{\otimes r} \otimes \mathbf{1}_{[(2k-2)/n,(2k/n)]}^{\otimes r} \rangle_{\mathfrak{H}^{\otimes 2}} |$$

$$= O(n^{-r}).$$

Consequently, when $r \neq 3$, we deduce

$$|R_n^{(r)}| \leq C n^{-(r+3)/2} \sum_{j,k=1}^{\lfloor n/2 \rfloor} |\rho(2j-2k)| = O(n^{-(r+1)/2}) \underset{n \to +\infty}{\longrightarrow} 0,$$

while when $r = 3$, we deduce

$$|R_n^{(3)}| \leq C n^{-1} \underset{n \to +\infty}{\longrightarrow} 0.$$

The proof of (4.4) is complete while the proof of (4.5) follows the same lines. Hence the proof of (4.2) is complete. □

LEMMA 4.3. *Let $f:\mathbb{R} \to \mathbb{R}$ be a function verifying* (H$_4$). *Set*

$$F_n = \sum_{k=1}^{\lfloor n/2 \rfloor} f(B_{(2k-1)/n})[(\Delta B_{(2k-1)/n})^2 - (\Delta B_{(2k-2)/n})^2].$$

*Then*

(4.6) $$F_n \xrightarrow[n\to\infty]{\text{stably}} \kappa \int_0^1 f(B_s)\,dW_s$$

*with $\kappa$ defined by (4.1) and where $W$ denotes a standard Brownian motion independent of $B$. Here, the stable convergence (4.6) is understood in the following sense: for any real number $\lambda$ and any $\sigma\{B\}$-measurable and integrable random variable $\xi$, we have that*

$$E(e^{i\lambda F_n}\xi) \xrightarrow[n\to\infty]{} E(e^{-\lambda^2\kappa^2/2 \int_0^1 f^2(B_s)\,ds}\xi).$$

PROOF. Since we follow exactly the proof of Theorem 3.2, we only describe the main ideas. First, observe that

$$F_n = \sum_{k=1}^{\lfloor n/2 \rfloor} f(B_{(2k-1)/n})(I_2(\delta_{(2k-1)/n}^{\otimes 2}) - I_2(\delta_{(2k-2)/n}^{\otimes 2})).$$

Here the analogue of Lemma 3.1 is

(4.7) $$\sup_{t\in[0,1]} \sum_{k=1}^{\lfloor n/2 \rfloor} |\langle \varepsilon_t, \delta_{(2k-1)/n}\rangle_{\mathfrak{H}}| \underset{n\to\infty}{=} O(1),$$

$$\sup_{t\in[0,1]} \sum_{k=1}^{\lfloor n/2 \rfloor} |\langle \varepsilon_t, \delta_{(2k-2)/n}\rangle_{\mathfrak{H}}| \underset{n\to\infty}{=} O(1),$$

(4.8) $$\left|\langle \varepsilon_{(2k-1)/n}, \delta_{(2k-1)/n}\rangle_{\mathfrak{H}}^2 - \frac{1}{4n}\right| \le \frac{\sqrt{2k} - \sqrt{2k-1}}{8n}$$

*and*

(4.9) $$\left|\langle \varepsilon_{(2k-1)/n}, \delta_{(2k-2)/n}\rangle_{\mathfrak{H}}^2 - \frac{1}{4n}\right| \le \frac{\sqrt{2k-1} - \sqrt{2k-2}}{4n}.$$

In fact, the bounds (4.7) are obtained by following the arguments presented in the proof of Lemma 3.1. The only difference is that in order to bound sums of the type $\sum_{k=1}^{\lfloor n/2 \rfloor} \sqrt{2k} - \sqrt{2k-1}$ (which are no more telescopic), we use

$$\sum_{k=1}^{\lfloor n/2 \rfloor} \sqrt{2k} - \sqrt{2k-1} \le \sum_{k=1}^{\lfloor n/2 \rfloor} \sqrt{2k} - \sqrt{2k-2} = \sqrt{2\lfloor n/2 \rfloor} \le \sqrt{n}.$$

As in Step 1 of the proof of Theorem 3.2, here we also have that $(F_n)$ is bounded in $L^2$. Consequently the sequence $(F_n, (B_t)_{t\in[0,1]})$ is tight in $\mathbb{R} \times \mathscr{C}([0,1])$. Assume that $(F_\infty, (B_t)_{t\in[0,1]})$ denotes the limit in law of a certain subsequence of $(F_n, (B_t)_{t\in[0,1]})$, denoted again by $(F_n, (B_t)_{t\in[0,1]})$. We have to prove that

$$(4.10) \qquad E(e^{i\lambda F_\infty} | (B_t)_{t\in[0,1]}) = \exp\left\{-\frac{\lambda^2}{2}\kappa^2 \int_0^1 f^2(B_s)\,ds\right\}.$$

We proceed as in Step 2 of the proof of Theorem 3.2. That is, (4.10) will be obtained by showing that for every random variable $\xi$ of the form (2.4) and every real number $\lambda$, we have

$$\lim_{n\to\infty} \phi'_n(\lambda) = -\lambda\kappa^2 E\left(e^{i\lambda F_\infty}\xi \int_0^1 f^2(B_s)\,ds\right),$$

where

$$\phi'_n(\lambda) := \frac{d}{d\lambda} E(e^{i\lambda F_n}\xi) = iE(F_n e^{i\lambda F_n}\xi), \qquad n \geq 1.$$

By the duality formula (2.5), we have that

$$\phi'_n(\lambda) = \sum_{k=1}^{\lfloor n/2 \rfloor} E(\langle D^2(f(B_{(2k-1)/n})e^{i\lambda F_n}\xi), \delta^{\otimes 2}_{(2k-1)/n} - \delta^{\otimes 2}_{(2k-2)/n}\rangle_{\mathfrak{H}^{\otimes 2}}).$$

The analogue of (3.8) is here:

$$\langle E(D^2(f(B_{(2k-1)/n})e^{i\lambda F_n}\xi)), \delta^{\otimes 2}_{(2k-1)/n} - \delta^{\otimes 2}_{(2k-2)/n}\rangle_{\mathfrak{H}^{\otimes 2}}$$
$$= E(f''(B_{(2k-1)/n})e^{i\lambda F_n}\xi)[\langle\varepsilon_{(2k-1)/n},\delta_{(2k-1)/n}\rangle_{\mathfrak{H}}^2 - \langle\varepsilon_{(2k-1)/n},\delta_{(2k-2)/n}\rangle_{\mathfrak{H}}^2]$$
$$+ 2i\lambda E(f'(B_{(2k-1)/n})e^{i\lambda F_n}\xi\langle DF_n,\delta_{(2k-1)/n}\rangle_{\mathfrak{H}})\langle\varepsilon_{(2k-1)/n},\delta_{(2k-1)/n}\rangle_{\mathfrak{H}}$$
$$- 2i\lambda E(f'(B_{(2k-1)/n})e^{i\lambda F_n}\xi\langle DF_n,\delta_{(2k-2)/n}\rangle_{\mathfrak{H}})\langle\varepsilon_{(2k-1)/n},\delta_{(2k-2)/n}\rangle_{\mathfrak{H}}$$
$$+ 2E(f'(B_{(2k-1)/n})e^{i\lambda F_n}\langle D\xi,\delta_{(2k-1)/n}\rangle_{\mathfrak{H}})\langle\varepsilon_{(2k-1)/n},\delta_{(2k-1)/n}\rangle_{\mathfrak{H}}$$
$$- 2E(f'(B_{(2k-1)/n})e^{i\lambda F_n}\langle D\xi,\delta_{(2k-2)/n}\rangle_{\mathfrak{H}})\langle\varepsilon_{(2k-1)/n},\delta_{(2k-2)/n}\rangle_{\mathfrak{H}}$$
$$- \lambda^2 E(f(B_{(2k-1)/n})e^{i\lambda F_n}\xi\langle DF_n,\delta_{(2k-1)/n}\rangle_{\mathfrak{H}}^2)$$
$$+ \lambda^2 E(f(B_{(2k-1)/n})e^{i\lambda F_n}\xi\langle DF_n,\delta_{(2k-2)/n}\rangle_{\mathfrak{H}}^2)$$
$$+ 2i\lambda E(f(B_{(2k-1)/n})e^{i\lambda F_n}\langle D\xi,\delta_{(2k-1)/n}\rangle_{\mathfrak{H}}\langle DF_n,\delta_{(2k-1)/n}\rangle_{\mathfrak{H}})$$
$$- 2i\lambda E(f(B_{(2k-1)/n})e^{i\lambda F_n}\langle D\xi,\delta_{(2k-2)/n}\rangle_{\mathfrak{H}}\langle DF_n,\delta_{(2k-2)/n}\rangle_{\mathfrak{H}})$$
$$+ i\lambda E(f(B_{(2k-1)/n})e^{i\lambda F_n}\xi\langle D^2 F_n,\delta^{\otimes 2}_{(2k-1)/n} - \delta^{\otimes 2}_{(2k-2)/n}\rangle_{\mathfrak{H}^{\otimes 2}})$$
$$+ E(f(B_{(2k-1)/n})e^{i\lambda F_n}\langle D^2\xi,\delta^{\otimes 2}_{(2k-1)/n} - \delta^{\otimes 2}_{(2k-2)/n}\rangle_{\mathfrak{H}^{\otimes 2}}).$$

As a consequence,

$$\phi'_n(\lambda) = -2\lambda \sum_{k,l=0}^{n-1} E(f(B_{(2k-1)/n})f(B_{(2l-1)/n})e^{i\lambda F_n}\xi)$$

(4.11)
$$\times \langle \delta^{\otimes 2}_{(2l-1)/n} - \delta^{\otimes 2}_{(2l-2)/n}, \delta^{\otimes 2}_{(2k-1)/n} - \delta^{\otimes 2}_{(2k-2)/n} \rangle_{\mathfrak{H}^{\otimes 2}}$$

$$+ i \sum_{k=0}^{n-1} r_{k,n},$$

where $r_{k,n}$ is given by

$$r_{k,n} = E[f''(B_{(2k-1)/n})e^{i\lambda F_n}\xi]$$

$$\times [\langle \varepsilon_{(2k-1)/n}, \delta_{(2k-1)/n} \rangle^2_{\mathfrak{H}} - \langle \varepsilon_{(2k-1)/n}, \delta_{(2k-2)/n} \rangle^2_{\mathfrak{H}}]$$

$$+ 2i\lambda E(f'(B_{(2k-1)/n})e^{i\lambda F_n}\xi \langle DF_n, \delta_{(2k-1)/n} \rangle_{\mathfrak{H}}) \langle \varepsilon_{(2k-1)/n}, \delta_{(2k-1)/n} \rangle_{\mathfrak{H}}$$

(4.12)
$$- 2i\lambda E(f'(B_{(2k-1)/n})e^{i\lambda F_n}\xi \langle DF_n, \delta_{(2k-2)/n} \rangle_{\mathfrak{H}}) \langle \varepsilon_{(2k-1)/n}, \delta_{(2k-2)/n} \rangle_{\mathfrak{H}}$$

$$+ 2E(f'(B_{(2k-1)/n})e^{i\lambda F_n} \langle D\xi, \delta_{(2k-1)/n} \rangle_{\mathfrak{H}}) \langle \varepsilon_{(2k-1)/n}, \delta_{(2k-1)/n} \rangle_{\mathfrak{H}}$$

$$- 2E(f'(B_{(2k-1)/n})e^{i\lambda F_n} \langle D\xi, \delta_{(2k-2)/n} \rangle_{\mathfrak{H}}) \langle \varepsilon_{(2k-1)/n}, \delta_{(2k-2)/n} \rangle_{\mathfrak{H}}$$

$$- \lambda^2 E(f(B_{(2k-1)/n})e^{i\lambda F_n}\xi \langle DF_n, \delta_{(2k-1)/n} \rangle^2_{\mathfrak{H}})$$

$$+ \lambda^2 E(f(B_{(2k-1)/n})e^{i\lambda F_n}\xi \langle DF_n, \delta_{(2k-2)/n} \rangle^2_{\mathfrak{H}})$$

$$+ 2i\lambda E(f(B_{(2k-1)/n})e^{i\lambda F_n} \langle D\xi, \delta_{(2k-1)/n} \rangle_{\mathfrak{H}} \langle DF_n, \delta_{(2k-1)/n} \rangle_{\mathfrak{H}})$$

$$- 2i\lambda E(f(B_{(2k-1)/n})e^{i\lambda F_n} \langle D\xi, \delta_{(2k-2)/n} \rangle_{\mathfrak{H}} \langle DF_n, \delta_{(2k-2)/n} \rangle_{\mathfrak{H}})$$

$$+ E(f(B_{(2k-1)/n})e^{i\lambda F_n} \langle D^2\xi, \delta^{\otimes 2}_{(2k-1)/n} - \delta^{\otimes 2}_{(2k-2)/n} \rangle_{\mathfrak{H}^{\otimes 2}})$$

$$+ i\lambda \sum_{l=1}^{\lfloor n/2 \rfloor} E(f(B_{(2k-1)/n})e^{i\lambda F_n}\xi f''(B_{(2l-1)/n})$$

$$\times (I_2(\delta^{\otimes 2}_{(2l-1)/n}) - I_2(\delta^{\otimes 2}_{(2l-2)/n})))$$

$$\times \langle \varepsilon^{\otimes 2}_{(2l-1)/n}, \delta^{\otimes 2}_{(2k-1)/n} - \delta^{\otimes 2}_{(2k-2)/n} \rangle_{\mathfrak{H}^{\otimes 2}}$$

$$+ 4i\lambda \sum_{l=1}^{\lfloor n/2 \rfloor} E(f(B_{(2k-1)/n})e^{i\lambda F_n}\xi f'(B_{(2l-1)/n})\Delta B_{(2l-1)/n})$$

$$\times \langle \varepsilon_{(2l-1)/n} \widetilde{\otimes} \delta_{(2l-1)/n}, \delta^{\otimes 2}_{(2k-1)/n} - \delta^{\otimes 2}_{(2k-2)/n} \rangle_{\mathfrak{H}^{\otimes 2}}$$

$$- 4i\lambda \sum_{l=1}^{\lfloor n/2 \rfloor} E(f(B_{(2k-1)/n})e^{i\lambda F_n}\xi f'(B_{(2l-1)/n})\Delta B_{(2l-2)/n})$$

$$\times \langle \varepsilon_{(2l-2)/n} \widetilde{\otimes} \delta_{(2l-2)/n}, \delta^{\otimes 2}_{(2k-1)/n} - \delta^{\otimes 2}_{(2k-2)/n} \rangle_{\mathfrak{H}^{\otimes 2}}$$

$$= \sum_{j=1}^{13} R^j_{k,n}.$$

The only difference with respect to (3.12) is that, this time, the term

$$i \sum_{k=0}^{n-1} E[f''(B_{(2k-1)/n})e^{i\lambda F_n}\xi][\langle \varepsilon_{(2k-1)/n}, \delta_{(2k-1)/n} \rangle^2_{\mathfrak{H}} - \langle \varepsilon_{(2k-1)/n}, \delta_{(2k-2)/n} \rangle^2_{\mathfrak{H}}]$$

corresponding to (3.15) is negligible. Indeed, we can write

$$\sum_{k=0}^{n-1} E[f''(B_{(2k-1)/n})e^{i\lambda F_n}\xi]$$

$$\times [\langle \varepsilon_{(2k-1)/n}, \delta_{(2k-1)/n} \rangle^2_{\mathfrak{H}} - \langle \varepsilon_{(2k-1)/n}, \delta_{(2k-2)/n} \rangle^2_{\mathfrak{H}}]$$

$$= \sum_{k=0}^{n-1} E(f''(B_{(2k-1)/n})e^{i\lambda F_n}\xi) \left[ \langle \varepsilon_{(2k-1)/n}, \delta_{(2k-1)/n} \rangle^2_{\mathfrak{H}} - \frac{1}{4n} \right]$$

$$- \sum_{k=0}^{n-1} E(f''(B_{(2k-1)/n})e^{i\lambda F_n}\xi) \left[ \langle \varepsilon_{(2k-1)/n}, \delta_{(2k-2)/n} \rangle^2_{\mathfrak{H}} - \frac{1}{4n} \right]$$

$$\xrightarrow[n\to\infty]{} 0 \qquad \text{by (4.8)–(4.9), under (H}_4\text{)}.$$

Moreover, exactly as in the proof of Theorem 3.2, one can show that

$$\lim_{n\to\infty} \sum_{k=1}^{\lfloor n/2 \rfloor} r_{k,n} = 0.$$

Consequently, we have

$$\lim_{n\to\infty} \phi'_n(\lambda)$$

$$= -2\lambda \lim_{n\to\infty} \sum_{k,l=1}^{\lfloor n/2 \rfloor} E(f(B_{(2k-1)/n})f(B_{(2l-1)/n})e^{i\lambda F_n}\xi)$$

$$\times \langle \delta^{\otimes 2}_{(2l-1)/n} - \delta^{\otimes 2}_{(2l-2)/n}, \delta^{\otimes 2}_{(2k-1)/n} - \delta^{\otimes 2}_{(2k-2)/n} \rangle_{\mathfrak{H}^{\otimes 2}}$$

$$= -\frac{\lambda}{2} \lim_{n\to\infty} \frac{1}{n} \sum_{k,l=1}^{\lfloor n/2 \rfloor} E(f(B_{(2k-1)/n})f(B_{(2l-1)/n})e^{i\lambda F_n}\xi)$$

$$\times (2\rho^2(2k-2l) - \rho^2(2l-2k+1) - \rho^2(2l-2k-1))$$

$$= -\frac{\lambda}{4} \sum_{r=-\infty}^{\infty} (2\rho^2(2r) - \rho^2(2r+1) - \rho^2(2r-1))$$

$$\times \lim_{n \to \infty} \frac{2}{n} \sum_{k=1 \vee (1-r)}^{\lfloor n/2 \rfloor \wedge (\lfloor n/2 \rfloor - r)} E(f(B_{(2k-1)/n})f(B_{(2k-1-2r)/n})e^{i\lambda F_n}\xi)$$

$$= -\lambda \kappa^2 \int_0^1 E(f^2(B_s)e^{i\lambda F_\infty}\xi)\,ds,$$

where $\kappa$ is defined by (4.1). In other words, (4.10) is shown and the proof of Lemma 4.3 is done. $\square$

**Acknowledgments.** Some of our computations are inspired by the first draft of [16]. We are grateful to David Nualart for letting us use them freely.

LABORATOIRE DE PROBABILITÉS
ET MODÈLES ALÉATOIRES
UNIVERSITÉ PARIS VI
BOÎTE COURRIER 188
4 PLACE JUSSIEU
75252 PARIS CEDEX 05
FRANCE
E-MAIL: inourdin@gmail.com

INSTITUT FÜR MATHEMATIK
HUMBOLDT-UNIVERSITÄT ZU BERLIN
UNTER DEN LINDEN 6
10099 BERLIN
GERMANY
E-MAIL: areveill@mathematik.hu-berlin.de